\documentclass[10pt]{article} 
\font\smallit=cmti10 
\font\smalltt=cmtt10

\usepackage{amssymb,latexsym,amsmath,epsfig,amsthm} 

\makeatletter

\renewcommand\section{\@startsection {section}{1}{\z@} 
{-30pt \@plus -1ex \@minus -.2ex}
{2.3ex \@plus.2ex} 
{\normalfont\normalsize\bfseries\boldmath}}

\renewcommand\subsection{\@startsection{subsection}{2}{\z@} 
{-3.25ex\@plus -1ex \@minus -.2ex}
{1.5ex \@plus .2ex} 
{\normalfont\normalsize\bfseries\boldmath}}

\renewcommand{\@seccntformat}[1]{\csname the#1\endcsname. }

\makeatother

\newtheorem{proposition}{Proposition}

\theoremstyle{definition} 
\newtheorem{definition}{Definition} 
\newtheorem{conjecture}{Conjecture} \newtheorem{remark}{Remark}

\begin{document}
\begin{center}
\uppercase{\bf Finite field models of 
polynomials 
interpolating Fourier 
coefficients of modular functions
for Hecke groups}
\vskip 20pt
{\bf Barry Brent}\\

{\tt barrybrent\@iphouse.com}\\
\end{center}
\vskip 20pt

\centerline{\smallit Received: , Revised: , Accepted: , Published: } 

\centerline{\bf Abstract}

\noindent
Following work of Raleigh and Akiyama
 ([31], [1]),
in our article [8]
we considered (among other objects)
families of weight zero 
meromorphic modular forms
$J_m$
for Hecke groups
$G(\lambda_m)$.
We conjectured 
that, for a certain
uniformizing variable $X_m$,
the $J_m$ have Fourier expansions
$J_m = 1/X_m + 
\sum_{n = 0}^{\infty} m^{-2n-2} A_n(m) X_m^n$,
where the $A_n(x)$ are polynomials
in $\mathbb{Q}[x]$. 
The present article is concerned with
models $\mathcal{A}_n[p](x)$ of the $A_n(x)$:
polynomials representing
self-maps of finite  fields
with characteristic $p$.
The main content is a conjecture
specifying $\mathcal{A}_n[p](x)$ up to
a multiplicative constant for certain
families of $n$ and $p$,
based on numerical experiments.

\pagestyle{myheadings}
\markright{\smalltt INTEGERS: 23 (2023)\hfill} 
\thispagestyle{empty}
\baselineskip=12.875pt
\vskip 30pt

\section{Introduction}
\subsection{Why Hecke groups?}
Here we describe
an old puzzle to
explain our interest in 
modular forms for Hecke
groups
and to advertise the
puzzle itself, which we still
have not solved, and indeed will not
address in the present
article.
\footnote{We have a draft
of a paper on this question 
in the folder 
``current draft'' in
our GitHub depository [5].}
C. L. Siegel ([37],[38])
established bounds on
the least positive integer
represented by a 
positive-definite
even unimodular quadratic form
in $2h$ variables
by first bounding
the exponent of the first
non-vanishing Fourier coefficient
for a level one entire 
classical modular form
$T_h$ of weight
$h$ such that the constant term
of $T_h$ is non-vanishing.
While working on
an extension of Siegel's result on 
the non-vanishing of the $T_h$ constant terms
to level two modular forms,
we came across the regularities described
in equations (1) and (2) in numerical
experiments. 
Let 
$\Delta$ denote the weight twelve
modular form for $SL(2,\mathbb{Z})$
that generates Ramanujan's tau function,
let $j$ be the usual Hauptmodul
normalized to have constant
term $744$, let
$d_b(n)$ be the sum of the 
digits in the base $b$
expansion of $n$, and
let $C(f)$ stand for the constant
term of the Fourier series
of $f$ in whatever uniformizing variable
happens to be in question. Then (apparently)
\begin{equation}
\text{ord}_2(C(j^k)) =\text{ord}_2(C(1/\Delta^k))
= 3d_2(k)
\end{equation}
and
\begin{equation}
\text{ord}_3(C(j^k)) = 
\text{ord}_3(C(1/\Delta^k)) =
d_3(k).
\end{equation}
We listed many functions displaying
analogous behavior in our 1998 article
[9],
together with a finer taxonomy based on 
congruences. (For example,
if $\text{ord}_3(n) = 1$, then $n$
may be congruent to $3$ or $6$
modulo $9$.)
Clearly,
if only we had proofs of these statements and
their analogues, we would
know that the constant terms of 
$1/\Delta^k, j^k$, and their analogues
are non-zero.
In the level one case Siegel
used different arguments
to establish the non-vanishing of
$T_h$ constant terms,
as we eventually did 
for their level two analogues
([9],[10].)
\newline \newline \noindent
We wondered about
the special role
of the primes $p = 2$ and $3$ in equations
(1) and (2): why these primes but not others?
(This is the 
aforementioned puzzle.) We looked for
patterns in the
$p$-orders of constant terms of 
$j$ and other modular forms for
$SL(2,\mathbb{Z})$ for $p$ 
larger than three.
Our search within 
$SL(2,\mathbb{Z})$ came up
empty, so 
we searched among the Hecke groups
$G(\lambda_n), n = 3, 4, ...$
for the following reasons.
The matrix
group $SL(2,\mathbb{Z})$ coincides 
with the Hecke group
$G(\lambda_3)$, discussed below. 
It is
isomorphic to the product of
cyclic groups $C_2*C_3$;
while in general  
 $G(\lambda_m) \cong C_2*C_m$ 
for $m = 3, 4, ....$
At first we hoped that primes
$p$ larger than three might
manifest behavior analogous
to that of three in equations (1) and (2)
in Hecke groups
$G(\lambda_p)$, but so far we have only found
some more complicated patterns.
\subsection{Existence of the interpolating 
polynomials.}
For $m = 3, 4, ...$, let 
$\lambda_m = 2 \cos \pi/m$ and
$J_m$ be a certain 
meromorphic modular 
form built from a particular
triangle function $\phi_m$
for the Hecke group $G(\lambda_m)$
with Fourier expansions
$
J_m(\tau) = \sum_{n=-1}^\infty a_n(m)q_m^n,
$
where
$q_m(\tau) = \exp 2 \pi i \tau/\lambda_m$.
The groups $G(\lambda_3)$ and 
$SL(2,\mathbb{Z})$ coincide.
\footnote{For further details, the reader
is referred to the books by 
 Carath{\'e}odory
[14], [15]
and by Berndt and Knopp
[4],
the articles of Lehner and Raleigh
[26], [31],
to the dissertation of Leo 
[27], and to a
summary, including pertinent references
to that material, in the
2021 article [8].  
Finally, the article by Hardy and Ramanujan
on expansions of
modular functions
is reprinted
in Ramanujan's \it Collected Papers \rm
[20].}
(There is more on triangle
functions in the next section,
and details of their construction are in 
[8].)
For  $n = -1, 0, 1, 2$ and $3$,
Raleigh [31]
gave polynomials $P_n(x)$ such that 
$a_{-1}(m)^n q_m^{2n+2} a_n(m) = P_n(m)$
for $m = 3, 4, ...,$
and conjectured that similar 
relations hold for all positive integers
$n$.
This was proved by Akiyama [1].
\newline \newline \noindent
In his 2021 article [8]
the present writer suggested that
such interpolating polynomials
for higher weight Hecke-group modular
forms should exist as well, and (acting
on his own suggestion) tentatively
identified some
of them by Lagrange interpolation. 
In this section, we
offer a more detailed existence argument.
Beyond this introduction,
the present article
does not require such an argument because we
will not construct models of 
interpolating polynomials
other than (more or less)
those of Raleigh and Akiyama. 
But a referee, whom we thank, asked us
to comment on future directions of this
investigation, and if the studies presented
here are to be extended to higher-weight
modular forms, 
the existence argument will be 
relevant. Our argument will
apply to Hecke-group modular forms.
The classical forms for $SL(2,\mathbb{Z})$
lie within that class.
\subsubsection{Hecke's theory of modular forms.}
Using the weight-raising properties of 
differentiation and the $J_m$, E.
Hecke constructed
certain families 
$\mathcal{H}$
comprising modular forms 
of positive weight
for each $G(\lambda_m)$
sharing certain properties 
([21], [4].)
(The weight of $g$ is not necessarily
constant within such a family.)
It seems apparent that Akiyama's 
result can be extended:
there should exist polynomials 
$Q_{\mathcal{H},n}(x)$ 
interpolating the coefficient of $X_m^n$
in the Fourier expansions of the
members of Hecke families
$\mathcal{H}$. 
\newline \newline \noindent
To make this precise, we review 
results of Hecke described in
the book of Berndt and Knopp 
[4]. 
By Theorem 3.1 in that book,
the region $B(\lambda_m)$ defined below is
a fundamental region for $G(\lambda_m)$.
\begin{definition}
\begin{enumerate}
\item Let $\tau_{\lambda_m}$ be the intersection
of the circle $|\tau| = 1$
with the line $\Re(\tau) = -\lambda_m/2$.
\item Let $B(\lambda_m)= \{\tau \in \mathbb{H}:
\Re(\tau) < \lambda_m/2,  |\tau| > 1\}$. 
\item Let $g_m(\tau)$ be the unique 
function guaranteed
to exist by the Riemann mapping function
mapping $B(\lambda_m)$ conformally
and one-to-one onto the upper half plane
such that $g_m$ takes  $\tau_{\lambda_m}$
to zero, $i$ to $1$, and $i\infty$ to itself.
(Berndt and Knopp, pages 47--48.) 
\item Let
 $$
 f_{\lambda_m}(\tau)
 :=\left \{
 \frac {g_m'(\tau)^2}
 {g_m(\tau)(g_m(\tau) -1)}
 \right \}^{1/(m-2)},
 $$
 $$
 f_{i,m}(\tau):=
 \left \{
 \frac {g_m'(\tau)^m}
 {g_m(\tau)^{m-1} (g_m(\tau) - 1)}
 \right \}^{1/(m-2)},
 $$
 and
 $$
 f_{\infty,m}(\tau) := 
 \left \{ \frac{g_m'(\tau)^{2m}}
 {g_m(\tau)^{2m-2}(g_m(\tau)-1)^m}
\right \}^{1/(m-2)}.
 $$
 \end{enumerate}
 \end{definition}
 \noindent
By Theorem 5.5 in Berndt and Knopp 
[4],
we know that the functions $f_{\lambda_m},
f_{i,m}$, and $f_{\infty,m}$ are modular
for $G(\lambda_m)$ with weights
$4/(m-2)$, $2m/(m-2)$, and $4m/(m-2)$,
respectively. 
(There is a
subtlety about the multiplier in
the functional equation for the
modularity of $f_{i,m}$ which we will pass over.)
\newline \newline \noindent
Because of its uniqueness,
we know that $g_m = J_m$ 
from equation 2 in Raleigh's article.
Therefore, corresponding to the three $f$'s,
we have the following definitions.
\begin{definition}
\begin{enumerate}
\item Let
$H_{\lambda, m}(\tau) :=$ 
 $$\left \{
 \frac {J_m'(\tau)^2}
 {J_m(\tau)(J_m(\tau) -1)}
 \right \}^{1/(m-2)}.$$
 \item Let $H_{\lambda,4,m} (\tau) :=
 H_{\lambda, m}(\tau)^{m-2}.$
 \end{enumerate}
 \end{definition}
 \begin{definition}
 \begin{enumerate}
 \item Let
  $H_{i,m}(\tau):=$ 
  $$
\left \{
 \frac {J_m'(\tau)^m}
 {J_m(\tau)^{m-1} (J_m(\tau) - 1)}
 \right \}^{1/(m-2)}.$$
 \item  Let $H_{i,6,m}(\tau):=$
  $$
\left \{
 \frac {J_m'(\tau)^m}
 {J_m(\tau)^{m-1} (J_m(\tau) - 1)}
 \right \}^{3/m}.$$
 \end{enumerate}
 \end{definition} 
 \begin{definition}
 \begin{enumerate}
 \item Let
$\Delta_{\infty,m}(\tau):=$
 $$
 \left \{ \frac{J_m'(\tau)^{2m}}
 {J_m(\tau)^{2m-2}(J_m(\tau)-1)^m}
\right \}^{1/(m-2)}.
 $$
 \item Let $\Delta_{\infty,12,m}(\tau):=$
 $$
 \left \{ \frac{J_m'(\tau)^{2m}}
 {J_m(\tau)^{2m-2}(J_m(\tau)-1)^m}
\right \}^{3/m}.
 $$
\item Let $\Delta^{\diamond}_m(\tau):=
H_{\lambda,m}(\tau)^3/J_m(\tau)$.
\item Let
$\Delta^{\diamond}_{12,m}(\tau):= 
H_{\lambda,4,m}^3(\tau)/J_m(\tau)$.
\item Let $\Delta^{\dagger}_m(\tau):=
H_{\lambda,4,m}(\tau)^3-H_{i,6,m}(\tau)^2$.
 \end{enumerate}
 \end{definition} 
 \begin{remark}
 It is
 easy to see from the definitions
 (for example, in [36]) that 
 in the classical case (subgoups of 
 $SL(2,\mathbb{Z})$), if
 $f$ and $g$ are modular for a particular group
 with weights $\omega_f$ and $\omega_g$,
 and $a$  is a rational number,
 then $fg$ and  $f^a$ are modular for the same group,
 with weights $\omega_f + \omega_g$ 
 and $a \cdot \omega_f$, respectively.
 These statements hold in the case of the Hecke groups
 as well.
 Therefore it follows from Berndt and Knopp's 
 Theorem 5.5 that 
  we have the following tables of weights:
 \begin{center}
\begin{tabular}{|c|c|c|c|} \hline
$H_{\lambda,m}$& $H_{\lambda,4,m}$ &$H_{i,m}$
& $H_{i,6,m}$  \\ \hline
 $4/(m-2)$&  $4$  & $2m/(m-2)$ & $6$ \\ \hline
\end{tabular} 
\end{center}
and
\begin{center}
\begin{tabular}{|c|c|c|c|c|} \hline
 $\Delta^{\diamond}_m$ 
&  $\Delta^{\diamond}_{12,m}$ & $\Delta_{\infty,m}$ 
& $\Delta_{\infty,12,m}$
& $\Delta^{\dagger}_m$\\ \hline
 $12/(m-2)$ &  $12$ &  
 $4m/(m-2)$ & $12$ & $12$\\ \hline
\end{tabular} 
\end{center}
 \end{remark}
\subsubsection{The existence argument.}
The main goal of the present article
is to describe conjectures about finite 
field models of the $J_m$
based on \it SageMath \rm 
experiments [33].
Let \footnote{Relevant files in [7]
are (1) the \it SageMath \rm Jupyter notebook
in which we generated Fourier expansions
of the $J_m$ and of the  $A_n(x)$, namely 
``capital-J make data file1jun21.ipnyb'';
(2) the notebook
``capital-J polynomials make file.ipynb''
in which we generated the data file
``run14jun21no14.txt'',
which is called in turn by many of our 
other notebooks;
(3) a table for $a_n(m)$
divided into several files:
``run2jun21no11'', ``run2jun21no12'', 
``run2jun21no13'' and ``run2jun21no14'';
(4) a \it Mathematica \rm notebook 
``conjecture 1.nb'' documenting 
the table's calculation;
a table for the $A_n(x)$ 
made in the same notebook 
 under file name 
``run20apr21no5'';
and (5) files
generated for
the same purposes 
in the SageMath Jupyter notebook
``conjecture1no1.ipnyb''
(which is also there.)}
$$
J_m(\tau) = 
\sum_{n = -1}^{\infty}a_n(m) q_m(\tau)^n.
$$
For integers $m \geq 3$ 
and $n = 0, 1, 2, 3$, 
$m^{-2n-2} a_{-1}(m)^{-n} A_n(m) = a_n(m)$,
where $A_n(x)$ is a polynomial with rational
coefficients and degree $2n+2$,
such that the coefficient of $x^n$
is zero when $n$ is odd 
(Raleigh, [31]).
As we said in the introduction,
similar relations
exist among the $a_n$ for 
all positive $n$.
\newline \newline \noindent
In section four of [8]
we wrote the Fourier expansions of 
the $J_m$ by replacing $q_m$ with 
another variable $X_m(\tau)$
 in the expansions. (The
 number $\tau$ is
a generic element of the upper half-plane.)
By Akiyama's theorem, we have a
series of the form
$\mathcal{J}(x,X_m) := 
\sum_{n=-1}^{\infty} \tilde{P}_n(x)X_m^n$
for polynomials $\tilde{P}_n(x)$ 
in $\mathbb{Q}[x]$
with the property that 
$J_m = \mathcal{J}(m,X_m)$.
A reader willing to take
that construction for granted can, 
for present purposes,
regard the expansion of $J_m$ in $X_m$ 
(truncated to $n$ terms) as the object defined 
in our \it SageMath \rm code
in the ``dictionary'' at the top of 
each notebook
as \tt{J(n,m)}\rm. 
\noindent Now, just to
provide the reader with some context, 
we normalize 
the $J_m$ themselves to obtain 
functions $j_m$ such that $j_3$
is (apparently) the usual  $j$ 
function.\footnote{Except to repeat
a conjecture from our 2021 article,
we do not
study the $j_m$ in this one.
But the reader will notice
in conjecture 1 below that 
that interpolating polynomials
behave a little more nicely for
$j_m$ than for $J_m$.}
We begin by defining an operator on
infinite series in $X_m$.
It has the effect when $m=3$
of recovering the Fourier series
of a variety of standard 
modular forms.\footnote{The 
substitution involved appears in 
[27]; see our article
[8] for a fuller 
acknowledgement of
our debt to Leo.}
\begin{definition}
Let $f = \sum_{n=a}^{\infty} k_n X_m^n$, 
where $k_n$ is a rational number 
for $n = a, a+1, ...,$
and $k_a \neq 0$. 
Let 
$g = \sum_{n=a}^{\infty} k_n (2^6 m^3 X_m)^n =
\sum_{n=a}^{\infty} \tilde{k}_n X_m^n$ (say). 
Then
$$\overline{f} := g/\tilde{k}_a.$$
\end{definition} 
\begin{definition} 
With $J_m(\tau) = 
1/X_m+\sum_{n = 0}^{\infty}a_n(m) X_m^n$,
we set\footnote{Some code for 
 $j_m$ Fourier expansions
appearing in
\it SageMath \rm notebooks
cited below was generated in
[7], notebook ``j from scratch.ipynb'', which
employs a ``dictionary''
(the definitions at the top of the notebook)
distinct from the corresponding dictionaries
in the notebooks where it is
reproduced.} 
$$j_m(\tau): = \overline{J_m} =
1/X_m +
\sum_{n\geq 0} c_m(n) X_m^n  \text{ (say).}
$$
\end{definition} 
\noindent
The Fourier expansion of $j_3$
is 
\footnote{
See equation (23) of 
Serre's book 
[36], 
section 3, and the \it SageMath \rm notebook 
``jpower constant term NewmanShanks 26oct22.ipynb''
in our repository [5].}
$$j_3(\tau) = 
1/X_3(\tau) + 744 + 196884 X_3(\tau) + 
21493760 X_3(\tau)^2 + ...,
$$
which matches the  standard expansion $j(\tau) = $
$$
1/\exp(2\pi i \tau) + 744 + 
196884 \exp(2\pi i \cdot  \tau) + 
21493760 \exp(2\pi i\cdot  2 \cdot \tau ) + ....
$$
\noindent
 We make the following 
 \begin{definition}
 Let $\mathcal{F} = \{f_3, ..., f_m, ...\}$
 where $f_m$ is modular for $G(\lambda_m)$.
 Then we write the Fourier expansion
 of $f_m^k$ in powers of $X_m$ as
$$
f_m(\tau)^k = \sum_n A_{\mathcal{F},k,m}(n)X_m^n.
$$
 \end{definition} \noindent
\begin{proposition}
Let 
$\mathcal{K} =
\{J_3, , J_4,...\}$
and $\overline{\mathcal{K}} =
\{j_3, j_4,....\}$
Then there exist polynomials
$Q_{\mathcal{K},k,n}(x)$ and
$Q_{\overline{\mathcal{K}},k,n}(x)$
in $\mathbb{Q}[x]$
such that 
$$J_m(\tau)^k = 
\sum_{n = -k}^{\infty} 
Q_{\mathcal{K},k,n}(m) X_m(\tau)^n$$ and
$$j_m(\tau)^k = 
\sum_{n = -k}^{\infty} 
Q_{\overline{\mathcal{K}},k,n}(m) 
X_m(\tau)^n.$$
In other words,
$A_{\mathcal{K},k,m}(n) = Q_{\mathcal{K},k,n}(m)$
and 
$A_{\overline{\mathcal{K}},k,m}(n) = 
Q_{\overline{\mathcal{K}},k,n}(m)$
for $k = 1, 2, ..., m = 3, 4, ...$,
and $n = -k, 1-k, ....$
\end{proposition} 
\noindent
For $k$ equal to one, 
the first claim is
just Akiyama's theorem and
the claim for $k$ not equal to
one is then obvious.
The second statement 
follows immediately.
\begin{proposition}
With $k$ as in proposition 1,
let $$\mathcal{H}
= \{H_{\lambda,m}\},
\{H_{\lambda,4,m}\},
\{H_{i,m}\},
\{H_{i,6,m}\},$$
$$
\{\Delta^{\diamond}_m\},
\{
\Delta^{\diamond}_{12,m}\},
\{\Delta_{\infty,m}\},
\text{or }
\{
\Delta^{\dagger}_m\},
$$
permitting $m$ to range over the integers 
greater than two.
Then there exist polynomials 
$Q_{\mathcal{H},k,n}(x)$
in $\mathbb{Q}[x]$
such that the elements $f_3, f_4, ...$ of 
$\mathcal{H}_k$
have Fourier expansions 
$$f_m(\tau) = 
\sum_n Q_{\mathcal{H},k,n}(m) X_m(\tau)^n.$$
\end{proposition} \noindent
For $k$ equal to one,
we justify this as follows.
After substituting $\mathcal{J}(x,X_m)$
(the series defined in the paragraph succeeding
remark 1 above)
for $J_m$ in the various clauses of 
definitions 2 - 4,
the right sides become rational 
functions of fractional powers
of various series in
powers of $X_m$ with coefficients
in $\mathbb{Q}[x]$,
which by purely formal operations
should be expressible as other
series in powers of $X_m$ with
coefficients in $\mathbb{Q}[x]$,
from which we 
recover Fourier expansions of 
each of the defined functions
by setting $x$ equal to $m$.
The statement for $k$ other
than one follows easily.
\subsection{Finite field models.}
\subsubsection{How do finite field models help?}
We want to understand as much as we can about 
the $A_n(x)$. The most
obvious ways to 
analyze polynomials are through their coefficients 
and their roots. 
We have not found patterns among the coefficients
in these polynomials. 
In experiments that
we have not published, 
on the other hand,
we did see that 
in some situations the roots of
relevant polynomials
appear to be confined
to the real 
axis, but we had no proofs.\footnote{
In some circumstances this is a significant
property. For example,
the proposition that 
all of the roots of Jensen polynomials
are real is equivalent to the
Riemann hypothesis ([30],
[19], [18], [17].)
}
We usually could not do better
than approximations.
The roots of polynomials
in a finite field, on the other
hand, can be enumerated
exactly, not just in theory but
(when the field is small
enough) 
in practice, simply by 
doing a brute-force search.
We found that (by clearing denominators
of a prime $p$) we arrived
at models of the $A_n(x)$
in fields of characteristic $p$
that displayed regularities
determined by the residue
class of $n$ modulo $p$.
These regularities
suggest, to the author,
at least,
that the $A_n$
in his 2021 article
were correctly
identified.  
They are also
evidence that the Fourier
coefficient of the original $J_m$
at $X^n$ is 
governed somehow by $n$ modulo $p$
for each prime $p$ and all $n$.
This
vague claim is the
moral of our story.\footnote{The
congruences for the coefficients $c(n)$ 
of the $j$ function discussed
in the articles 
[2], [24], and [25],
and also in in Serre's book 
[36] (chapter VII, section 3)
are clearly related to this claim.}
\subsubsection{A sketch of
the models.}
If $P(x)$ is in $\mathbb{Q}[x]$,
we write  $\mathcal{K}_p(P(x))$
for a certain polynomial 
 in $\mathbb{Q}[x]$
that agrees with $P(x)$ 
 up to a multiplicative constant.
The definition of the map $\mathcal{K}_p$ 
guarantees that $\text{ord}_p$
is non-negative on all coefficients 
of $\mathcal{K}_p(P(x))$--a property that 
allows us to define its models
in finite fields of characteristic $p$.
We write
$\mathbb{A}(n,p)$ for the smallest 
splitting field (in a weak sense)
of $A_n(x)$
over $\mathbb{F}_p$. We also write
 $\mathcal{A}_n[p](x)$
for  $\mathcal{K}_p(A_n(x))$
considered 
as a polynomial self-map of $\mathbb{A}(n,p)$.
These are the models of our title.
We study
the $\mathcal{A}_n[p](x)$
with numerical experiments
and make conjectures describing their
roots completely
for primes $p$ less than or equal to seven. 
\newline \newline \noindent
To recover the $A_n(m)$
from the $\mathcal{A}_n[p](x)$,
our experiments suggest that
it would be most useful to have
in hand $\mathcal{A}_n[p](x)$
for all primes $p$ less than
or equal to the smallest prime
greater than $n$.
But the feasibility of computing with $p$
declines with the size of $\mathbb{A}(n,p)$,
which becomes large 
for $n$ in certain residue classes modulo $p$.
Consequently, our conjectures only address
primes  less than or equal to seven. If
our conjectures are valid and the behavior
of the $\mathbb{A}(n,p)$ is similar for all 
larger primes as well, then the Fourier
coefficients $A_n(m)$ 
of the $J_m = 1/X_m + 
\sum_{n = 0}^{\infty} A_n(m) X_m^n$ behave in
a uniform way that is independent of
the Hecke group $G(\lambda_m)$ and depends only on 
congruences satisfied by the indices $n$. 
\newline \newline \noindent
Our identification of their roots for 
$p$ less than or equal to seven only specifies the
$\mathcal{A}_n[p](x)$ themselves up 
to a multiplicative constant, and we have not been 
able to understand how these constants vary with $n$,
even for fixed small values of $p$. On the other hand,
the original polynomials $A_n(x)$ behaved much better
in this regard in our experiments. We were
able to write them as the product of
a product of monic polynomials and an
explicitly known rational number. (See conjectures
1 and 2 below.) This seems
to lessen the urgency of solving the
problem of understanding
more fully the multiplicative constants associated 
with the $\mathcal{A}_n[p](x)$.
\section{Triangle functions.}
The material in the present section 
was sketched more fully in 
[8], which also
includes citations
to an even more  thorough exposition of much of it 
in the second volume of Carath\'eodory [15].
Let $\mathbb{Z}, \mathbb{Q}, 
\mathbb{C}$ and $\mathbb{H}$ 
denote, respectively, the set of rational integers,
the set of rational
numbers, the set  of 
complex numbers, and the set of
complex numbers
with positive imaginary parts. 
(We will reserve 
the letter $\tau$ for elements of the
upper half-plane, and $z$
for generic complex numbers.)
We write
$\mathbb{H}^* = \mathbb{H} \cup \mathbb{Q} \cup 
\{i \infty\}$, and we equip $\mathbb{H}^*$
with the Poincar{\'e} metric. Figures
$T$ made by three geodesics
of $\mathbb{H}^*$
are called hyperbolic or
circular-arc triangles.
Let $\lambda_m = 2 \cos \pi/m$.
For $m = 3, 4, ...$,
we define the Hecke group
$G(\lambda_m)$
as the discrete group generated by 
the maps $z \rightarrow -1/z$
and $ z \rightarrow z + \lambda_m$.
The full modular group $SL(2,\mathbb{Z})$
is identical to $G(\lambda_3)$.
\newline  \newline \noindent
For our purposes,
Schwarz triangles $T$ are 
hyperbolic triangles in $\mathbb{H}^*$
with certain restrictions
on the angles at the vertices.
From a Euclidean point of view, 
their sides are
vertical rays, 
segments of vertical rays, 
semicircles 
orthogonal to the real axis
and meeting it
at points $(r,0)$
with $r$ rational,
or arcs of such semicircles.
We choose
$\lambda, \mu$ and $\nu$,
all non-negative, such that
$\lambda + \mu + \nu <1$;
then the angles of $T$
are $\lambda \pi, \mu \pi$,
and $\nu \pi$.
By reflecting $T$ 
across one of its edges, we
get another Schwarz triangle. 
The reflection between
two triangles in $\mathbb{H^*}$
is effected by 
a M\"obius transformation,
so the orbit of $T$
under repeated reflections
is associated to a 
 collection of 
 M\"obius transformations.
The group 
 generated by 
these transformations is a
triangle group.\footnote{In the terminology of 
Isant and Grau [3], page 45,
these are Fuchsian groups of 
signature $(0;1/\lambda,1/\mu,1/\nu)$.
We have not as yet examined the
possibility of extending the experiments of 
our article [8] (and of the 
present article) to a broader class
of Fuchsian groups. Movasati [28]
 provides
relevant code for such studies,
based on results of the article by Doran, Gannon
and himself [16].}
By the Riemann Mapping Theorem
there is a conformal, onto 
map
$\phi: T \mapsto \mathbb{H}^*$
called a triangle function.
\newline \newline \noindent
Hecke groups 
are triangle groups 
$H$ that act 
properly discontinuously 
on $\mathbb{H}$ [21].
This means that for compact 
$K \subset \mathbb{H}$, the set
$\{\mu \in H$ s.t. \hskip -.05in
$K \cap \mu(K) \neq \emptyset\}$
is finite.
Recall that 
$G(\lambda_m)$ is the Hecke group
generated by the maps 
$z \mapsto -1/z$
and $z \mapsto z + \lambda_m$.
Hecke established 
that $G(\lambda_m)$ has the structure
of a free product 
of cyclic groups
$C_2 * C_m$,
generalizing the relation
([36], [13])
$SL(2,\mathbb{Z}) = C_2 * C_3$.
\newline \newline \noindent
Let $\rho = -\exp(-\pi i/m) = 
-\cos(\pi/m) + i\sin(\pi/m)$, and
let $T_m \subset \mathbb{H}^*$ 
denote the 
hyperbolic triangle with vertices
$\rho, i$, and $i\infty$. 
The corresponding angles are 
$\pi/m, \pi/2$ and $0$
respectively. Let 
$\phi_{\lambda_m}$ be a 
triangle function 
for $T_m$. The
function $\phi_{\lambda_m}$ 
has a pole at $i\infty$
and period $\lambda_m$.
For $P, Q \in \mathbb{H}^*$,
let us us write $P \equiv_H Q$ when 
$\mu \in H$ and $Q = \mu(P)$.
Then $\phi_{\lambda_m}$ extends to 
a function
$J_m: \mathbb{H}^* \rightarrow \mathbb{H}^*$ by
declaring that
$J_m(P) = J_m(Q)$
if and only if 
$P \equiv_H Q$.
$J_m$ is a modular function
(a meromorphic modular form of weight zero)
for $G(\lambda_m)$.
\newline \newline \noindent
Schwarz ([34] or [35]),
Lehner [26]
and Raleigh [31]
studied Schwarz triangle functions, 
which map  hyperbolic triangles $T$ in the 
 extended upper half $z$-plane  onto
 the extended upper half 
 $w$-plane.
 For certain $T = T_m$, 
a triangle function 
$\phi_{\lambda_m}: T \to \mathbb{H}^*$
 extends to a map
$J_m: \mathbb{H}^* \rightarrow \mathbb{H}^*$
 invariant under
modular transformations
from $G(\lambda_m)$.
\section{Polynomial interpolation of 
Fourier coefficients of
modular functions for 
Hecke groups.}
\noindent
We reproduce parts of conjecture 1 and
conjecture 2 from our earlier
article [8].
\begin{conjecture}
Let 
the Fourier expansion of 
$j_m(\tau)$ be
$$j_m = 1/X_m +
\sum_{n\geq 0} c_m(n) X_m^n.$$
\rm  
\begin{enumerate}
\item 
For each integer $n$
greater than $-2$,
there exists a polynomial
$C_n(x) \in \mathbb{Q}[x]$
that satisfies the relation
$c_m(n) = C_n(m)$ for $m = 3, 4, ...$
\footnote{[7], notebook
``conjecture 1.nb''.}
\item 
Let $\{\phi_{-1}, \phi_0, ...\} = \{1, 24, ...\}$
be the McKay-Thompson series of class 4A 
[40]. For some
 degree $2n$, irreducible, 
 monic polynomial $\gamma_n(x)$
in $\mathbb{Q}[x]$ we have \footnote{
[7], notebook ``conjecture 1 clause 2.ipynb'}:
$$
C_n(x) = 
\phi_n \cdot (x-2)(x+2)x^{n+1} \gamma_n(x).
$$
\item 
The function $j_3$ is identical to the
modular function on $SL(2,\mathbb{Z})$
usually denoted $j$.\footnote{[7],
notebook ``conjecture 1 cause 3.nb''.}
\end{enumerate}
\end{conjecture} 
\begin{conjecture}
Let \footnote{Relevant documents in [7]
are notebooks
``conjecture 2.nb'',
``conjecture2no1.ipynb'',
``capital-J make data file1jun21.ipynb''
and associated data files.} 
the Fourier expansion of $J_m(\tau)$
be
$$
J_m = 
\sum_{n = -1}^{\infty}a_m(n) X_m^n.$$
\begin{enumerate}
\item \footnote{For clause 1, see 
    [7], notebooks ``conjecture 2.nb'',
    ``conjecture 2 clause 1b.ipynb'', and 
``conjecture 2 clause 1b no2.ipynb''.}
    \begin{enumerate}
        \item There exist polynomials $A_n(x)$
    such that
    $A_{-1}(x) \equiv 1$,
    $A_0(x) = 3x^2+4$, 
    $A_1(x) = 69x^4 - 8x^2 - 48$,
    and $A_n(m) = m^{2n+2}a_m(n)$ 
    for $m = 3, 4, ....$.
    \item 
Let $C_n(x)$ be as in Conjecture 1. Then
$$A_n(x) = 2^{-6n-6} x^{-n-1} C_n(x).$$
\end{enumerate}
    \item  
Let $\pi_n$ be the
    set of prime numbers
    dividing the denominator
    of at least one non-zero
    coefficient of $A_n$.
    Then the following statements
    are true \footnote{[7],
notebook ``conjecture 2 clause 2 
w code 14jun21.ipynb''.}
    \begin{enumerate}
    \item $\pi_2 = \{3\}$.
      \item If $\pi_n$ is ordered by size, it contains
        no gaps. That is, if $p$ and $p'$
        are consecutive elements of  $\pi_n$ 
        with $p = p_k$ and $p' = p_j$,
        then $j = k + 1$.
            \item If $n$ is an odd prime, 
            then
            $$\pi_n = \{2, ...,k,...,p\}_{k \text{ prime}}$$ where
            $p$ is the greatest prime less than $n$.
            \item If $n$ is composite and
            $n+1$ is prime, then 
            $$\pi_n = \{2...,k, ...,n+1\}_{k \text{ prime}}.$$
            \item If $n$ and $n+1$ are both composite, 
            then  
            $$\pi_n = \{2...,k, ...,p\}_{k \text{ prime}}$$
            where $p$ is the greatest prime less than $n$.
    \end{enumerate}
\end{enumerate}
\end{conjecture}
\begin{remark}
    Comparing the notation of clause 1(a) of  conjecture
   2 with that of proposition 1, we see that
    $A_n(m) = m^{2 n + 2}Q_{\mathcal{K},1,n}(m)
    =  m^{2 n + 2}A_{\mathcal{K},1,m}(n)
    $.
\end{remark}
\section{Finite-field models of interpolating polynomials.} 
Interpreting a polynomial $f(x)$ 
as a self-map on
a finite field $\mathbb{F}_{p^k}$
is not possible
if the prime $p$ divides the denominator of
one of the coefficients of $f(x)$. 
It appears that
all primes less than
$n$ occur in the denominators
of coefficients of the
interpolating polynomials $A_n(x)$ 
and that, as a result,
we cannot make such an interpretation
when $p$ is less than $n$.
To get around this difficulty, we
clear denominators  to obtain
polynomials, either with integer coefficients,
or at least with denominators not divisible 
by a specified prime number,
before applying 
maps from $\mathbb{Q}[x]$
to $\mathbb{F}_{p^k}[x]$
defined below to arrive at finite field models of
the interpolating polynomials.\footnote{In 
practice, these maps are
\it SageMath \rm coercions. Coercion is a 
concept from type theory. For a 
recent discussion, see Buzzard's
preprint [11].
For \it SageMath\rm's own explanations
of its coercion routines, see
the documents [42] and [22].
\it Caveat: \rm In the article,
we define the maps 
 mathematically; coercion
in \it SageMath\rm, on the other hand,
is defined by code.
Our conjectures
are based on the validity 
(which is left to the judgement 
of the reader)
of the correspondence
between these definitions.}
\begin{definition}
\begin{enumerate}
\item 
If $F$ is a finite field of characteristic $p$, 
let $0_F$ 
and $1_F$ denote
the additive and multiplicative identities
of $F$, respectively.
Furthermore if $n$ is a positive integer,
let $n_F$ 
denote the sum in $F$
of $n$ copies of $1_F$,
and let $(-n)_F$ be the
additive inverse in $F$ of $n_F$.
If $n_1$ and $n_2$ are positive integers
such that $\text{ord}_p(n_2) \geq 0$
let $(n_1/n_2)_F$ be the quotient
in $F$ of $(n_1)_F$ and $(n_2)_F$.
We refer to the 
map $x \mapsto x_F$
as coercion. 
\item 
\begin{enumerate}
\item Let
$\mathbb{A}(n,p)$
be the smallest 
characteristic $p$ field in which $A_n(x)$
splits into factors of degree zero or one.
\item Let $\mathbb{A}(n,p)^*$ be the (cyclic)
multiplicative group 
$\mathbb{A}(n,p) \backslash \{0_F\}$.
\end{enumerate}
\item
\begin{enumerate}
\item Let
$s_A(n,p)$ denote the vector space
dimension $[\mathbb{A}(n,p):\mathbb{F}_p]$.
\item Integers $n_1$ and $n_2$
are $(A,p)$-equivalent if
$n_1 \equiv n_2$ modulo $p$
and \newline $s_A(n_1,p)  = s_A(n_2,p)$.
\end{enumerate}
\item \begin{enumerate}
    \item If
$P(x) = \sum_{n \in I} (N_n/D_n) x^n \in \mathbb{Q}[x]$
where $N_n$ and $D_n$ are relatively prime
integers, $d(P(x)) := \text{lcm} \{D_n\}_{n \in I}$.
\item With $\mu_p(P(x)) = p^{\text{ord}_p(d(P(x))}$,
$$K(P(x)):= d(P(x)) \times P(x)$$
and 
$$\mathcal{K}_p(P(x)):= \mu_p(P(x)) \times P(x).$$
\end{enumerate}
\item  Here we define 
finite field models of the $A_n(x)$. 
\begin{enumerate}
    \item 
Let $\alpha_j, j =  0, 1, 2,...,u,$
be rational integers, and let 
$K(A_n(x)) = \sum_{j=0}^u \alpha_j x^j$.
Then 
$$
A_n[p](x) := \sum_{j=0}^u(\alpha_j)_{\hskip -.02in F} 
\thinspace x^j.
$$
where  $F = \mathbb{A}(n,p)$.
\item 
Let $\alpha_j^*, j =  0, 1, 2,...,u,$
be rational numbers  and let 
$\mathcal{K}_p(A_n(x)) = \sum_{j=0}^u \alpha_j^* x^j$.
By construction,
$\text{ord}_p(\alpha_j^*) \geq 0$.
Then 
$$
\mathcal{A}_n[p](x) := 
\sum_{j=0}^u(\alpha_j^*)_{\hskip -.02in F} \thinspace x^j.
$$
where  $F = \mathbb{A}(n,p)$.
\end{enumerate}
\item
\begin{enumerate}
\item Let $\tt{mod}$ \hskip -.05in $(n,p)$   
be the element $\rho$ of
 $\{0, 1, ..., p-1 \}$ such that \newline
$n \equiv \rho$ \hskip .01in (\hskip -.1in $\mod p)$.
\item We set
$\delta(n,p):= (n-\tt{mod}$ \hskip -.05in $(n,p))/p$.
\item We set
$a(n,p) :=$ the leading coefficient not divisible 
by $p$ in $K(A_n(x))$.
\item We set
$a^*(n,p) :=$ the leading coefficient not divisible 
by $p$ in $\mathcal{K}_p(A_n(x))$.
\item We set
$r(n,p) := \tt{mod}\it (a(n,p),p)$.
\item We set
$r^*(n,p) := \tt{mod}\it (a^*(n,p),p)$.
\item We set
$\alpha(n,p):=r(n,p)_F$ where $F = \mathbb{A}(n,p)$.
\item We set
$\alpha^*(n,p):=r^*(n,p)_F$ where $F = \mathbb{A}(n,p)$.
\item We write the Frobenius map as
$\phi_p: F \rightarrow F$ with
$\phi_p(s) := s^p.$
\item For $s$ in $F$, the $\phi_p$-orbit of $s$
is written $\mathcal{O}_p$.
\end{enumerate}
\end{enumerate}
\end{definition}
\noindent
\begin{remark}
\begin{enumerate}
\item The models $\mathcal{A}_n[p](x)$
and $A_n[p](x)$ agree up to  
multiplicative constants. Therefore, 
the assertions summarized in tables 1 - 3 below
are valid for both of them.
We require the $A_n[p](x)$
in order to make the argument 
in remark 5 concerning 
Frobenius orbits of the roots. 
But there is a loss of
information in $A_n[p](x)$
from the  clearing 
of denominators by the 
$K$-operator which can only be
reflected in the relevant 
multiplicative constant. 
To be more precise,
 $K$ takes distinct polynomials
 in $\mathbb{Q}[x]$ to a single
 element of $\mathbb{Z}[x]$,
and so it would seem to be 
difficult to constrain 
$A_n(x)$ using $K(A_n(x))$
or the corresponding models $A_n[p](x)$.
On the other hand, 
let the set of prime numbers
$\pi_n$ be as in
conjecture 2, clause 2 and
let $\nu_n$ be the largest element
of $\pi_n$. Then the obstacles to
deriving conditions on $A_n(x)$ from 
information about the images
$\{\mathcal{K}_p(A_n(x))\}_{p \leq \nu_n}$ and 
the corresponding models
$\{\mathcal{A}_n[p](x)\}_{p \leq \nu_n}$
would seem to be less formidable.
\item
We have not considered the question
whether or not the $\mathbb{A}(n,p)$
are splitting fields in the full sense of,
say, Definition A.5.7  in the book [29].
(The extra condition is this: for a splitting field $S$
 of a polynomial $P(x)$, the roots
of $P(x)$ generate $S$  over the base field.)
\item Typically 
when $n_1$ and $n_2$ are $(A,p)$-equivalent,
the polynomials $\mathcal{A}_{n_1}[p](x)$
and $\mathcal{A}_{n_2}[p](x)$ are distinct
(tables 1-3.)
Of course, the set of self-maps of a finite set
is also finite, so the infinite family
$\{\mathcal{A}_n[p](x)\}_n$
represents only a finite set of self-maps
on the $\mathbb{A}(n,p)$,
and there is a natural equivalence relation
on $\{\mathcal{A}_n[p](x)\}_n$ induced by
this fact; we will study it below.
\end{enumerate}
\end{remark}
\noindent
The conjectures below are based on
numerical experiments documented in 
our repository [6] in the range of
$n$-values that were accessible
to our tests,
namely, $-1 \leq n \leq 200$. 
\begin{conjecture}\footnote{[6],
notebooks with titles ``conjecture 4 ...''
 (\it sic.\rm) for various primes $p$; these notebooks
serve several purposes, hence the file names.
In the file names,
``res'' stands for ``residue'', 
\it i.e.\rm, \tt{mod}\rm$\hskip -.00in (n,p)$,
and
``sd'' stands for
``splitting degree'', \it i.e.\rm, $s_A(n,p)$. } 
Let $p_1$ and $p_2$ be prime numbers
such that $p_1 < p_2$.\footnote{[6],
notebook 'conjecture 3 for $A_n$.ipynb'.}
\begin{enumerate}
\item \footnote{[6], 
Jupyter notebook
``numerical term on extension fields 17may22.ipynb''} 
\begin{enumerate}
\item
If $A_n[p](x)$ is factored over $\mathbb{A}(n,p)$
as a field element $\gamma$ times a 
product of monic polynomials, 
then $\gamma = \alpha(n,p)$;
by construction, $\gamma$ belongs to 
$\{1_{\mathbb{A}(n,p)}, 2_{\mathbb{A}(n,p)}, ...,
(p-1)_{\mathbb{A}(n,p)}\}$.
\item If $\mathcal{A}_n[p](x)$ is 
factored over $\mathbb{A}(n,p)$
as a field element $\gamma$ times a 
product of monic polynomials, 
then $\gamma = \alpha^*(n,p)$ and
$\gamma$ belongs to \newline
$\{1_{\mathbb{A}(n,p)}, 2_{\mathbb{A}(n,p)}, ...,
(p-1)_{\mathbb{A}(n,p)}\}$.
\end{enumerate}
\item If $p$ is a prime greater than $3, p < n$
and  $p \leq 17$,
then $s_A(n,p)$ is determined by 
$\tt{mod}$ \hskip -.05in $(n,p)$;
moreover, $s_A(n,2) \equiv 1$ identically. 
Here are tables for
$p = 3, 5$ and $7$.\footnote{[6], 
notebook ``conjecture 3 clause 2.ipynb''}
\begin{center}
\begin{tabular}{ |c|c|c|c|c|c| } \hline
$\tt{mod}$ \hskip -.05in $(n,3)$ & $0$ & $1$ & $2$\\ \hline
 $s_A(n,3)$ & $1$ & $1, 2, 3$ & $1, 2, 3, 4$\\ \hline
\end{tabular}
\end{center}
\begin{center}
\begin{tabular}{ |c|c|c|c|c|c| } 
 \hline
$\tt{mod}$ \hskip -.05in $(n,5)$ & $0$ & $1$ & 
$2$ & $3$ & $4$ \\ \hline
 $s_A(n,5)$  & $1,2$ & $4$ & $4$ & $6$ & $1$  \\ \hline
\end{tabular}
\end{center}
\begin{center}
\begin{tabular}{ |c|c|c|c|c|c|c|c| } \hline
 $\tt{mod}$ \hskip -.05in $(n,7)$ & $0$ & $1$ & $2$ & $3$ 
 & $4$ & $5$ & $6$\\ \hline
 $s_A(n,7)$ & $4$ & $4$ & $2$ & $2$ & $6$ & $5$ & $1$  \\ \hline
\end{tabular}
\end{center}
\item If $p$ is a prime number greater than three and 
$n \equiv p-1$ \rm (\hskip .005in mod  $p$),
\it then  $s_A(n,p) = 1$.
\item 
\begin{enumerate}
\item Either $s_A(0,p) = 1$ or $s_A(0,p) = 2$.
\item If $p$ is a prime number greater than three and 
$s_A(0,p) = 1$, then \newline
$p \equiv 1$ \rm (\hskip .005in mod  $3$).\it
\item If $p$ is a prime number greater than three and 
$s_A(0,p) = 2$, then \newline
$p \equiv 2$ \rm (\hskip .005in mod  $3$).
\end{enumerate}
\end{enumerate}
\end{conjecture}
\begin{remark}
\begin{enumerate}
    \item We have not found statements 
    like clause 4 
about $s_A(n,p)$ for
values of $n$ other than zero. 
\item We found clause 4 after reading
Sloane's A002476 [39] and 
related pages.
\end{enumerate}
\end{remark}
\begin{remark}
By construction, the coefficients of 
$K(A_n(x))$  are
rational integers; hence
the coefficients of
$A_n[p](x)$ 
lie in the prime sub-field
of $\mathbb{A}(n,p)$.
Let $f(x)$ be a polynomial self-map
of $\mathbb{A}(n,p)$ with coefficients
 in its prime sub-field. 
If $t$ is a generator of $\mathbb{A}(n,p)^*$, then
$t^{kp}$ is a root of $f(x)$ whenever $t^k$ is 
such a root. (Proof: the Frobenius map 
$\phi_p: x\mapsto x^p$ 
is an automorphism of $\mathbb{A}(n,p)$
and fixes its prime sub-field.\footnote{For 
example, see Chapter VII, section 5 
of Lang's book [23].}
So applying $\phi_p$ to the equation
$f(x) = 0$  gives 
$f(\phi_p(x)) = 0$.) 
Consequently, the roots of $A_n[p](x)$ 
and $\mathcal{A}_n[p]$  are (identical)
collections of complete $\phi_p$ orbits.
(In particular, each root of the form 
$n_{\mathbb{A}(n,p)}$ 
is a $\phi_p$ fixed point.)
\end{remark}
\newpage
\begin{conjecture}\footnote{This conjecture is
based on  numerical data
available in our GitHub repository [6] in files
with filenames beginning ``conjecture 4 p = (*),
res=(**), sd=(***)'' for various integers (*), etc.
When these are loaded, they show titles beginning 
``conjecture 4 clause 4''.
The reader should disregard the phrase ``clause 4''.}
We state conjectures  on
$\mathcal{A}_n[p](x)$
for given $(A,p)$ equivalence class 
data and $p = 2, 3, 5$ and $7$. 
Let $t(n,p) = t$ (say)
be a generator of $\mathbb{A}(n,p)^*$. 
Excepting members $n_F$ of
the prime sub-field, we display 
the non-zero roots of 
$\mathcal{A}_n[p](x)$ as
 $\phi_p$-orbits  $\mathcal{O}_p(t^k)$
 for some $k$.\footnote{ A reader who
 consults our \it Sagemath \rm notebooks in
the repository [6] will find that we kept track
 of the action of Frobenius by
 computing the (appropriately padded)
 base $p$ expansions of 
 the discrete logarithms $k$;
 when $\phi_p(t^k) = t^j$,
 the expansion of $j$ is the image
 of the expansion of $k$ under a
 cyclic permutation.
 For example, in the case
 $p = 3, n \equiv 1 
 \text{\thinspace (mod}  3),
 s_A(n,3) = 2$ (row 4 of table 2),
 our \it SageMath \rm code
 outputs a root $t+2_F$ with
 discrete logarithm base $3$ 
 expansion $(1,2)$. This tells us
 that $t + 2 = t^7$ in
 $\mathbb{A}(n,3)$ because 
 $(1,2)$ is the base $3$ 
 expansion of $7$,
 and that the length of
 $\mathcal{O}_3(t+2)$ is $2$,
 because the length of the orbit of
 $(1,2)$ under cyclic 
 permutation is $2$.}
 We do not make a general
 conjecture regarding the sizes
 of the Frobenius orbits, but
for $p = 2, 3,$ or $5$,
let $r$ be a non-zero root of
$\mathcal{A}_n[p](x)$.
If $r$ belongs to the
prime sub-field of $\mathbb{A}(n,p)$
so that, for some rational
integer $n, r = n_F$,
then $r$ is a fixed point of $\phi_p$:
$\#\mathcal{O}_p(r) = 1$.
Otherwise,
$\#\mathcal{O}_p(r) = s_A(n,p)$.
\newline \newline \noindent
The situation for $p=7$ 
is more complicated.
If we write $r = t^k$, 
denote the base $p$ expansion of $k$
as $X_p(k)$, 
 padded with zeros,
necessary, so that its
length is equal to $s_A(n,p)$,
then for any prime $p$
the length of $\mathcal{O}_p(r)$ is 
the same as the length of of the
orbit under cyclic permutations 
of $X_p(k)$.
For the smaller $p$, $X_p(k)$ exhibits no
internal symmetries, but $X_7(k)$
does. For example,
if $n \equiv 4 (\text{mod \thinspace} 7$)
and $s_A(n,7) = 6$,
we find that 
$\mathcal{A}_n[7]$ has roots
$r_1=t^k$ for $k = 29412$ and
$r_2 = t^k$ for
$k = 88236$. 
Here $X_7(29412) = (5,1,5,1,5, 1)$
and $X_7(88236) = (1, 5, 1, 5, 1, 5)$.
These expansions comprise
a complete orbit under cyclic permutation.
Consequently, $r_1$ and $r_2$
comprise a complete orbit under $\phi_7$.
On the other hand, by a similar analysis,
$s_A(n,p) = 1$
implies that all roots
of $\mathcal{A}_n[p]$
are fixed points of $\phi_p$,
and $s_A(n,p) = 2$
implies that all roots $r$ of 
$\mathcal{A}_n[p]$ are either
fixed by $\phi_p$ or satisfy
$\# \mathcal{O}_p(r) = 2$.
Roots $r = t^k$ fixed by $\phi_p$
 must belong to 
the prime sub-field
(for example, the notes of Ca$\tilde{\rm{n}}$ez
[12], pp. 27-28).
\newline \newline \noindent
We have dropped 
constant factors from our tables.
For example, $3_F (x-1_F)$ 
would be listed simply as $x-1_F$.
\footnote{What we think 
might be true
about the constant factors
is stated in conjecture 3. 
We cannot specify them \it ab initio\rm.}
We use the $\delta(n,p)$
notation defined above and write 
$\delta(n,p) = \delta$; 
the arguments will be clear.
The letter ``u'' means unrestricted.
When we wish to abbreviate, ``$F$'' 
means $\mathbb{A}(n,p)$.
\begin{enumerate}
    \item 
 $\mathcal{A}_n[2] = x^{2n+2}.$
 \newpage
\item  Table for $p = 3$:
\begin{center}
\begin{tabular}{ |c|c|c|c| } \hline
$\tt{mod}$ \hskip -.05in $(n,3)$ 
 & $s_A(n,3)$

 &  $\mathcal{A}_n[3](x)$\\ \hline
$0$ 
& $1$

& 
   \begin{tabular}{@{}c@{}}
 \\
   $x^2(x-1_F)^{2\delta}(x-2_F)^{2\delta}$\\
  \\
 \end{tabular}
   \\ \hline
  $1$ 
& $1$ 

& 
   \begin{tabular}{@{}c@{}}
 \\
   $x^2(x-1_F)^{2\delta}(x-2_F)^{2\delta}$\\
  \\
 \end{tabular}
   \\ \hline
 $1$ 
& $1$ 

&
   \begin{tabular}{@{}c@{}}
 \\
   $x^2(x-1_F)^{2\delta-2}(x-2_F)^{2\delta-2}$\\
  \\
  \end{tabular}  \\ \hline
 $1$ 
& $2$

& \begin{tabular}{@{}c@{}}
\\
 $x^8(x-1_F)^{2\delta-6}
   (x-2_F)^{2\delta-6}$ \\
\\ $\times \prod_{s \in \mathcal{O}_3(t)}(x-s)^2$\\
\\ $\times \prod_{s \in \mathcal{O}_3(t^2)}(x-s)^2$\\
\\ $\times \prod_{s \in \mathcal{O}_3(t^7)}(x-s)^2$\\
  \\
 \end{tabular} \\\hline
  $2$ 
 & $1$ 
 
 & 
 \begin{tabular}{@{}c@{}}
 \\
 $(x-1_F)^{2\delta+1}(x-2_F)^{2\delta+1}$ \\
 \\
 \end{tabular} 
 
 \\ \hline
  $2$ 
  & $1$ 
 
 & \begin{tabular}{@{}c@{}}
 \\
 $x^6(x-1_F)^{2\delta-1}
   (x-2_F)^{2\delta-1}$ \\
   \\
 \end{tabular} \\\hline
 \end{tabular}
\newline \newline 
\sc{table 1: 
$\mathcal{A}_n[3]$ up to
constant factors}
\end{center}
\noindent
(In the ambiguous cases, we have not
deciphered the conditions that choose between
the polynomials listed in table 1 for 
$\mathcal{A}_n[3]$.)
 \newpage
\item Table for $p = 5$:
\begin{center}
\begin{tabular}{ |c|c|c|c|c| } \hline
$\tt{mod}$ \hskip -.05in $(n,5)$ 
 & $s_A(n,5)$ 
 &  $\mathcal{A}_n[5](x)$\\ \hline
 
 $0$ 
& $1$

 & \begin{tabular}{@{}c@{}}\\
 $(x-1_F)^{2\delta+2}
   (x-4_F)^{2\delta+2} \thinspace \times$\\
   \\
  $(x-2_F)^{2\delta}
   (x-3_F)^{2\delta}$\\
   \\
 \end{tabular}
   \\ \hline
$1$ 
& $4$

& \begin{tabular}{@{}c@{}}\\
 $(x-1_F)^{2\delta}
   (x-4_F)^{2\delta} \thinspace \times$\\
   \\
  $(x-2_F)^{2\delta}
   (x-3_F)^{2\delta} \times$\\
   \\
   $\prod_{s \in 
   \mathcal{O}_5(t^{91})}(x-s)$\\
   \\
 \end{tabular}  \\\hline
 $2$ 
& $4$

& \begin{tabular}{@{}c@{}}\\
 $(x-1_F)^{2\delta}
   (x-4_F)^{2\delta} \thinspace \times$\\
   \\
  $(x-2_F)^{2\delta+1}
   (x-3_F)^{2\delta+1} \times$\\
   \\
   $\prod_{s \in 
   \mathcal{O}_5(t^{169})}(x-s)$\\
   \\
 \end{tabular}  \\ \hline
 
$3$ 
& $6$

& \begin{tabular}{@{}c@{}}\\
 $(x-1_F)^{2\delta}
   (x-4_F)^{2\delta} \thinspace \times$\\
   \\
  $(x-2_F)^{2\delta+1}
   (x-3_F)^{2\delta+1} \times$\\
   \\
   $\prod_{s \in 
   \mathcal{O}_5(t^{2961})}(x-s)$\\
   \\
 \end{tabular}  \\\hline
 $4$ 
& $1$ 

& \begin{tabular}{@{}c@{}}\\
 $(x-1_F)^{2\delta+2}
   (x-2_F)^{2\delta+2} \thinspace \times$\\
   \\
  $(x-3_F)^{2\delta+2}
   (x-4_F)^{2\delta+2}$\\
   \\
 \end{tabular}  \\\hline
\end{tabular}
\newline \newline 
\sc{table 2: 
$\mathcal{A}_n[5]$ up to
constant factors.}
\end{center}
\newpage
\item 
Table for $p=7$:
\begin{center}
\begin{tabular}{ |c|c|c|c|c| } \hline
$\tt{mod}$ \hskip -.05in $(n,7)$ &  
 $s_A(n,7)$
 &  $\mathcal{A}_n[7](x)$\\ \hline
$0$ & $4$ 
   & \begin{tabular}{@{}c@{}}
   \\
 $(x-1_F)^{2\delta-2}
   (x-6_F)^{2\delta-2} \times$\\
   \\
    $(x-2_F)^{2\delta}
   (x-5_F)^{2\delta} \times$\\
   \\
   $(x-3_F)^{2\delta+1}
   (x-4_F)^{2\delta+1} \times$\\
   \\
   $\prod_{s \in 
   \mathcal{O}_7(t^{173})}(x-s)$\\
   \\ 
    $\prod_{s \in 
   \mathcal{O}_7(t^{260})}(x-s)\times$\\
   \\ 
 \end{tabular}
\\\hline
$1$ & $4$ 

   & \begin{tabular}{@{}c@{}}
   \\
 $(x-1_F)^{2\delta}
   (x-2_F)^{2\delta} \times$\\
   \\
    $(x-3_F)^{2\delta}
   (x-3_F)^{2\delta} \times$\\
   \\
   $(x-5_F)^{2\delta}
   (x-6_F)^{2\delta} \times$\\
   \\
   $\prod_{s \in 
   \mathcal{O}_7(t^{75})}(x-s)$\\
   \\ 
   \end{tabular} \\ \hline
   $2$ & $2$ 
   & \begin{tabular}{@{}c@{}}
   \\
 $(x-1_F)^{2\delta}
   (x-3_F)^{2\delta} \times$\\
   \\
    $(x-4_F)^{2\delta}
   (x-6_F)^{2\delta} \times$\\
   \\
   $(x-2_F)^{2\delta+1}
   (x-5_F)^{2\delta+1} \times$\\
   \\
   $\prod_{s \in 
   \mathcal{O}_7(t^4)}(x-s)^2$\\
   \\ 
 \end{tabular}
\\\hline
 \end{tabular}
\newline \newline
\sc{\noindent part one of table 3: 
$\mathcal{A}_n[7]$ up to
constant factors.}
\end{center}
\newpage
This is a continuation of the preceding table:
\begin{center} 
\begin{tabular}{ |c|c|c|c|c| } \hline
$\tt{mod}$ \hskip -.05in $(n,7)$ &  
 $s_A(n,7)$
 &  $\mathcal{A}_n[7](x)$\\ \hline
$3$ & $2$ 
   & \begin{tabular}{@{}c@{}}
   \\
 $(x-1_F)^{2\delta}
   (x-3_F)^{2\delta} \times$\\
   \\
    $(x-4_F)^{2\delta}
   (x-6_F)^{2\delta} \times$\\
   \\
   $(x-2_F)^{2\delta+1}
   (x-5_F)^{2\delta+1} \times$\\
   \\
   $\prod_{s \in 
   \mathcal{O}_7(t^7)}(x-s) \times$\\
   \\ 
   $\prod_{s \in 
   \mathcal{O}_7(t^{12})}(x-s) \times$\\
   \\  
   $\prod_{s \in 
   \mathcal{O}_7(t^{25})}(x-s)$\\
   \\  
 \end{tabular}\\ \hline
  $4$ & $6$ 
   & \begin{tabular}{@{}c@{}}
   \\
 $(x-1_F)^{2\delta}
   (x-3_F)^{2\delta} \times$\\
   \\
    $(x-4_F)^{2\delta}
   (x-6_F)^{2\delta} \times$\\
   \\
   $(x-2_F)^{2\delta+1}
   (x-5_F)^{2\delta+1} \times$\\
   \\
   $\prod_{s \in 
   \mathcal{O}_7(t^{29412})}(x-s) \times$\\
   \\
   $\prod_{s \in 
   \mathcal{O}_7(t^{41280})}(x-s) \times$\\
   \\ 
   $\prod_{s \in 
   \mathcal{O}_7(t^{81528})}(x-s)$\\
   \\ 
 \end{tabular} \\\hline
 \end{tabular}
 \vskip .2in
\sc{\noindent part two of table 3: 
$\mathcal{A}_n[7]$ up to
constant factors.}
\end{center}
\newpage
This is part three of the preceding table:
 \begin{center} 
\begin{tabular}{ |c|c|c|c|c| } \hline
$\tt{mod}$ \hskip -.05in $(n,7)$ &  
 $s_A(n,7)$
 &  $\mathcal{A}_n[7](x)$\\ \hline
$5$ & $5$ 
   & \begin{tabular}{@{}c@{}}
   \\
 $(x-1_F)^{2\delta}
   (x-3_F)^{2\delta} \times$\\
   \\
    $(x-4_F)^{2\delta}
   (x-6_F)^{2\delta} \times$\\
   \\
   $(x-2_F)^{2\delta+1}
   (x-5_F)^{2\delta+1} \times$\\
   \\
   $\prod_{s \in 
   \mathcal{O}_7(t^{1513})}(x-s) \times$\\
   \\ 
   $\prod_{s \in 
   \mathcal{O}_7(t^{11020})}(x-s)$\\
   \\ 
 \end{tabular}
\\\hline
$6$ 
& $1$ 
   & \begin{tabular}{@{}c@{}}\\
 $(x-1_F)^{2\delta+2}
   (x-2_F)^{2\delta+2} \times$\\
   \\
    $(x-3_F)^{2\delta+2}
   (x-4_F)^{2\delta+2} \times$\\
   \\
   $(x-5_F)^{2\delta+2}
   (x-6_F)^{2\delta+2}$\\
\\  
\end{tabular} \\ \hline
\end{tabular} 
\vskip .2in
\sc{\noindent part three of table 3: 
$\mathcal{A}_n[7]$ up to
constant factors.}
\end{center}
\noindent
In the case of the prime $p = 7$, 
there was only one example for which 
$n \equiv 0 \thinspace (\text{mod} \thinspace \thinspace 7)$ 
and $s_A(n,7) = 1$
in the range of our observations. Namely,
with $F = \mathbb{A}(0,7)$,
$\mathcal{A}_0[7] = 3_F (x-1_F)(x-6_F)$.
For the other cases, we have table 3 above.
By remarks and footnotes above, 
it is only necessary to
state the value of $\# \mathcal{O}_7$
on a root of  $\mathcal{A}_n[7]$
when $s_A(n,7)$ is larger than $2$.
For roots of $\mathcal{A}_n[7]$
where $n \equiv 4 \thinspace 
\text{\rm (mod \it} 7)$,
$\#\mathcal{O}_7(t^{29412}) = 2$,
$\#\mathcal{O}_7(t^{41290}) = 3$,
and
$\#\mathcal{O}_7(t^{81528}) = 3$.
For roots of $\mathcal{A}_n[7]$
where $n \equiv 5 \thinspace
\text{\rm (mod \it} 7)$,
$\#\mathcal{O}_7(t^{1513}) = 5$
and
$\#\mathcal{O}_7(t^{11020}) = 5$.
\end{enumerate}
\end{conjecture}



\begin{thebibliography}{1}\footnotesize

\bibitem{A} Shigeki Akiyama, A note on Hecke’s absolute invariants,  
{\it J. Ramanujan Math. Soc.} 
{\bf 7.1} \rm (1992), 65–81. 
\bibitem{AO} A. O. L. Atkin and J. N. O’Brien, 
Some properties of $p(n)$ and $c(n)$ modulo powers of 13,
{\it Trans. Amer. Math. Soc.} {\bf 126.3} \rm (1967), 442–459. 
\bibitem{IG} P. B. Isant and A. T. Grau, Uniformization of triangle modular curves, 
{\it Publicacions matema`tiques} \rm (2007), 0043–106. 
\bibitem{BK} B. C. Berndt and M. I. Knopp, 
{\it Hecke’s theory of modular forms and Dirichlet series}, \rm 
Vol. 5. World Scientific, 2008. 
\bibitem{BB1} B. Brent, constant terms repository, \newline
URL https://github.com/barry314159a/NewmanShanks. 
\bibitem{BB2} B. Brent, finite-models repository, \newline
URL https://github.com/barry314159a/trianglefunctions. 
\bibitem{BB3} B. Brent, interpolations repository, \newline URL 
https://github.com/barry314159a/interpolations. 
\bibitem{BB4} B. Brent, Polynomial interpolation of modular forms for 
Hecke groups, 
{\it Integers} {\bf 21} \rm (2021), \#A118.
URL: 
http://math.colgate.edu/(tilde)integers/v118/v118.pdf. 
\bibitem{BB5} B. Brent, Quadratic minima and modular forms, 
{\it Experiment. Math.} 
{\bf 7.3} \rm (1998), 257–274. 
\bibitem{BB6} B.Brent, Quadratic minima and modular forms II, 
{\it Acta Arith.} {\bf 96.4}\rm, (2001), 381–387.
\bibitem{B} K. Buzzard. What is the point of computers? 
A question for pure mathematicians, 
preprint, \tt{arXiv:2112.11598}. \rm
\bibitem{C1} S. Ca$\tilde{\rm{n}}$ez. Abstract Algebra 
Northwestern University Lecture Notes. 
\newline The following URL should be entered on 
a single line: \newline
https://sites.math.northwestern.edu/(tilde)scanez/
\newline
courses/331/notes/lecture-notes-331-3.pdf. 
\bibitem{C2} I. N. Ca\"ngul, The group structure of Hecke groups 
H($\lambda q$), 
{\it Turkish J. Math.} {\bf 20.2} \rm (1996), 203–207. 
\bibitem{C3} C. Carath\'eodory ({\it tr.} F. Steinhardt), 
{\it Theory of functions of a complex variable, 
Second English Edition,
Volume 1 }\rm   Chelsea Publishing Company, 1958. 
\bibitem{C4} C. Carath\'eodory ({\it tr.} F. Steinhardt), 
{\it Theory of functions of a complex variable, 
Second English Edition. Volume 2 }\rm, 
Chelsea Publishing Company, 1981.
\bibitem{D} C. F. Doran \it et al, \rm Automorphic 
forms for triangle groups, preprint, \newline
\tt{arXiv:1307.4372}. \rm
\bibitem{G1} M. Griffin \it et al, \rm Jensen Polynomials 
for the Riemann Xi Function,
preprint, \tt{arXiv:1910.01227}.  \rm
\bibitem{G2} M. Griffin \it et al, \rm Jensen polynomials 
for the Riemann xi-function,
{\it Adv. Math.} {\bf  397} \rm (2022), 108186. 
\bibitem{G3} M. Griffin \it et al.\rm, Jensen polynomials for 
the Riemann zeta function and 
other sequences,
{\it Proc. Natl. Acad. Sci. USA} {\bf 116.23} \rm (2019), 11103–11110. 
\bibitem{HR} G. H. Hardy and S. Ramanujan, 
On the coefficients in the expansions of certain modular functions, 
in {\it Collected papers of 
Srinivasa Ramanujan}, AMS Chelsea Publishing,
Providence, RI, 1962.
\newline \rm URL
https://royalsocietypublishing.org/doi/pdf/10.1098/rspa.1918.0056. 
\bibitem{H} E. Hecke, \"Uber die bestimmung dirichletscher 
reihen durch ihre funktionalgleichung, 
{\it Math. Ann.} {\bf 112.1} \rm (1936), 664–699. 
\bibitem{K} S. King. ``Sage’s category and coercion framework.'' 
\newline (The URL below should be entered as a single line.)
\newline
https://doc.sagemath.org/html/en/thematic(understroke)tutorials/
\newline
coercion(understroke)and(understroke)categories.html.
\bibitem{La} S. Lang. {\it Algebra} Springer Science and Business Media, 
Berlin, 2012. 
\bibitem{Le} J. Lehner, Divisibility properties of the 
Fourier coefficients of 
the modular \newline invariant 
$j(\tau)$, {\it Amer. J. Math.} {\bf 71.1} (1949), 136–148. 
\bibitem{Le2} J. Lehner, Further congruence properties of 
the Fourier coefficients 
of the \newline modular invariant $j(\tau)$, 
{\it Amer. J. Math.} {\bf 71.2}(1949), pp. 373–386. 
\bibitem{Le3} J. Lehner, Note on the Schwarz triangle functions, 
{\it Pacific J. Math.}
{\bf 4.2} (1954), 243–249. 
\bibitem{Leo} J. G. Leo, Fourier coefficients of 
triangle functions, Ph.D. thesis, 
\newline URL http://halfaya.org/ucla/research/thesis.pdf, 2008. 
\bibitem{M} H. Movasati, Modular forms/Automorphic forms for triangle groups, 
\newline 
URL https://w3.impa.br/
(tilde)hossein/computerprogramming.html. 
\bibitem{MM} G. L. Mullen and C. Mummert, 
{\it Finite Fields and Applications},
Stud. Math. Libr., Amer. Math. Soc., 2007. 
\bibitem{PJ} G. Polya and J. L. W. V. Jensen, {\it \"Uber die 
algebraisch-funktionentheoretische Untersuchungen von JLWV Jensen}, 
AF Høst, 1927. 
\bibitem{R1} J. Raleigh, On the Fourier coefficients of triangle 
functions, {\it Acta Arith.} 
{\bf 8} (1962), 107–111. 
\bibitem{R2} S. Ramanujan, {\it Collected papers of Srinivasa 
Ramanujan}, (Hardy, G. H., Seshu, P. V., and Wilson, 
B. M., {\it eds.}), Cambridge University Press, 2015, 310–321. 
\bibitem{SD} The Sage Developers, SageMath, the 
Sage Mathematics Software System (Version 9.6). 
\newline URL https://www.sagemath.org. 2021. 
\bibitem{Sc1} H. A. Schwarz, {\it Ueber diejenigen \"Falle, in welchen die 
Gaussische hypergeometrische Reihe eine algebraische Function 
ihres vierten Elementes darstellt}, Walter de Gruyter, 1873. 
\bibitem{Sc2} H. A. Schwarz, Ueber diejenigen \"Falle, in welchen die 
Gaussische hypergeometrische Reihe eine algebraische Function 
ihres vierten Elementes darstellt, 
{\it Journal f\"ur die reine und angewandte Mathe- matik}, 
{\bf 75} (1873), 292–335. 
\bibitem{Se} J.-P.Serre, {\it A course in arithmetic}, Springer-Verlag, 1970. 
\bibitem{Si1} C. L. Siegel, {\it Berechnung von Zetafunktionen an 
ganzzahligen Stellen}, 
Vandenhoeck and Ruprecht, 1969. 
\bibitem{SR} C. L. Siegel and K. G. Ramanathan, 
{\it Advanced analytic number theory}, 
Tata Institute of Fundamental Research Bombay, 1980, 249–268. 
\bibitem{S} N. J. A. Sloane, Primes of the form $6n+1$, 
{\it The On-Line Encyclopedia of Integer Sequences},
URL http://oeis.org/A002476. 
\bibitem{So} M. Somos, McKay-Thompson series of class 4A for 
the Monster group with $a(0) = 24$, 
{\it The On-Line Encyclopedia of Integer Sequences}, 
URL https://oeis.org/A097340.
\bibitem{SDT} The Sage Development Team. Parents, conversion and 
coercion. \newline URL
https://doc.sagemath.org/html/en/tutorial/tour
(understroke)coercion.html. 
\end{thebibliography}
\end{document}